\documentstyle{article}\begin{document}\centerline{\bf\Large Some results on Theory of Infinite Series and Divisor Sums}\vskip .2in

\centerline{\bf N.D. Bagis}
\centerline{\bf Aristotele University of Thessaloniki}
\centerline{\bf Thessaloniki Greece}
\centerline{\bf nikosbagis@hotmail.gr}

\[
\]

\centerline{\bf Abstract}
In this article we present certain formulas involving arithmetical functions. In the first part we study properties of sums and product formulas for general type of arithmetic functions. In the second part we apply these formulas to the study of Jacobi elliptic theta functions theory.

\section{General Theorems and Preparations}

\textbf{Proposition 1.}\\
If $x$ is positive real number and $f$ is analytic in $(-1,1)$ with $f(0)=0$, then
\begin{equation}
\exp\left(\int^{x}_{\infty}f(e^{-t})dt\right)=\prod^{\infty}_{n=1}(1-e^{-nx})^{\frac{1}{n}\sum_{d|n}\frac{f^{(d)}(0)}{d!}\mu\left(\frac{n}{d}\right)}, 
\end{equation}
where $x>0$ and $\mu$ is the Moebius-$\mu$ arithmetic function (see [1]). The function $\mu(n)$ take the values $(-1)^r$ when $n$ is square free and product of $r$ primes, else is $0$. Also $\mu(1)=1$.\\
\\
\textbf{Proof.}\\
Because $f(0)=0$ and $f$ analytic in $(-1,1)$, the integral $\int^{x}_{\infty}f(e^{-t})dt$ exists for every $x>0$. We assume that exists arithmetic function $X(n)$ such that: 
\begin{equation}
\exp\left(\int^{x}_{\infty}f(e^{-t})dt\right)=\prod^{\infty}_{n=1}(1-e^{-nx})^{X(n)}
\end{equation}
we will determine this function $X$.\\
Taking logarithms in both sides of (2) we have
$$
\int^{x}_{\infty}f(e^{-t})dt=\sum^{\infty}_{n=1}X(n)\log(1-e^{-nx})=-\sum^{\infty}_{n=1}X(n)\sum^{\infty}_{m=1}\frac{e^{-mx}}{m}=
$$
$$
=-\sum^{\infty}_{n,m=1}X(n)n\frac{e^{-mnx}}{mn}=-\sum^{\infty}_{n=1}\frac{e^{-nx}}{n}\sum_{d|n}X(d)d\eqno{:(A)}$$
Derivating (A) we get
$$
f(x)=\sum^{\infty}_{n=1}e^{-nx}\sum_{d|n}X(d)d\eqno{:(B)}
$$
But from analytic property of $f$ in $(-1,1)$ we have $$f(x)=\sum^{\infty}_{n=1}\frac{f^{(n)}(0)}{n!}x^n$$ and consequently
$$
f(e^{-x})=\sum^{\infty}_{n=1}\frac{f^{(n)}(0)}{n!}e^{-nx}
$$
Therefore from (B) and the above relation it must be 
$$
\frac{f^{(n)}(0)}{n!}=\sum_{d|n}X(d)d
$$
By applying the Moebius inversion theorem (see [1]) we get
$$
X(n)=\frac{1}{n}\sum_{d|n}\frac{f^{(d)}(0)}{d!}\mu\left(\frac{n}{d}\right)
$$
This completes the proof. 
Note also that holds the following similar expression
\begin{equation}
e^{-f(q)}=\prod^{\infty}_{n=1}\left(1-q^n\right)^{\frac{1}{n}\sum_{d|n}\frac{f^{(d)}(0)}{\Gamma(d)}\mu\left(\frac{n}{d}\right)}
\end{equation}
\\
\textbf{Examples on Proposition 1.}\\
\textbf{1)} If $f(x)=x$ then $\frac{f^{(n)}(0)}{n!}=\delta_n$, $n=1,2,3,\ldots$ i.e $\delta_1=1$, 0 else. Hence  $X(n)=\frac{1}{n}\sum_{d|n}\delta_d\mu\left(\frac{n}{d}\right)=\frac{\mu(n)}{n}$ and
\begin{equation}
\prod^{\infty}_{n=1}(1-q^n)^{\frac{\mu(n)}{n}}=e^{-q}
\end{equation}
\textbf{2)}
Let $\frac{f^{(n)}(0)}{n!}=n$, $n=1,2,3,\ldots$, then $f(x)=\frac{x}{(x-1)^2}$ and $X(n)=\frac{1}{n}\sum_{d|n}d\mu\left(\frac{n}{d}\right)=\frac{\phi(n)}{n}$. Where $\phi(n)=\sum_{d|n}d\mu\left(\frac{n}{d}\right)$ is Euler's phi arithmetic function.\\Hence
\begin{equation}
\prod^{\infty}_{n=1}(1-q^n)^{\frac{\phi(n)}{n}}=e^{\frac{q}{q-1}}
\end{equation}
\textbf{3)}
Let $\frac{f^{(n)}(0)}{n!}=\frac{\mu(n)}{n^{\nu}}$, $n=1,2,3,\ldots$, then $f(x)=\sum^{\infty}_{n=1}\frac{\mu(n)x^n}{n^{\nu}}$ and $X(n)=\frac{1}{n}\sum_{d|n}d^{-\nu}\mu(d)\mu\left(\frac{n}{d}\right)=\frac{\sigma^{(-1)}_{-\nu}(n)}{n}$, i.e $\sigma_{\nu}(n)=\sum_{d|n}d^{\nu}$ is the sum of the $\nu$-th power of divisors of $n$ and $\sigma^{(-1)}_{\nu}$ is its arithmetic inverse.\\ This means  $\sum_{d|n}\sigma_{\nu}(d)\sigma^{(-1)}_{\nu}\left(\frac{n}{d}\right)=\delta_n$. 
\begin{equation}
\prod^{\infty}_{n=1}(1-q^n)^{\frac{\sigma^{(-1)}_{-\nu}(n)}{n}}=\exp\left(-\sum^{\infty}_{n=1}\frac{\mu(n)q^n}{n^{\nu+1}}\right),
\end{equation}
where $\sigma^{(-1)}_{\nu}(n)$ is the aritrhmetic inverse of $\sigma_{\nu}(n)$.\\ 
\\
\textbf{Theorem 1.}\\
When $a,b>0$ and $f$ is analytic in $(-1,1)$ then
\begin{equation}
\prod^{\infty}_{n=1}\left(\frac{1-e^{-nb}}{1-e^{-na}}\right)^{\frac{1}{n}\sum_{d|n}\frac{f^{(d)}(0)}{d!}\mu\left(\frac{n}{d}\right)}=\exp\left(\int^{b}_{a}f(e^{-t})dt\right)
\end{equation}
\textbf{Proof.}\\
Easy consequence of Proposition 1.\\
\\
\textbf{Proposition 2.}\\
If $a$ is positive real number then
\begin{equation}
\sum^{\infty}_{n=1}\frac{\sum_{d|n}\frac{f^{(d)}(0)}{d!}\mu\left(\frac{n}{d}\right)}{e^{na}-1}=f\left(e^{-a}\right)
\end{equation}
\textbf{Proof.}\\ 
Set $x=a>0$ in (1) and take the logarithmic derivative in both sides with respect to $a$.\\
\\ 
\textbf{Proposition 3.}\\
\begin{equation}
\sum^{\infty}_{n=1}\frac{\sum_{d|n}\frac{f^{(d)}(0)}{d!}\mu\left(\frac{n}{d}\right)}{e^{na}+1}=-2f\left(e^{-2a}\right)+f\left(e^{-a}\right)
\end{equation}
\textbf{Proof.}\\
Set $x=a$ and $x=2a$ in (1) to take two relations, divide them. Take the logarithms and derivate. After a few simplifications we get (9).\\
\\  
\textbf{Proposition 4.}\\
If $A(n)$ is arbitrary arithmetic function, then for $x>0$ we have
\begin{equation}
\frac{d^{\nu}}{dx^{\nu}}\left(\sum^{\infty}_{n=1}\frac{\sum_{d|n}A(d)\mu\left(\frac{n}{d}\right)}{e^{nx}-1}\right)=\sum^{\infty}_{n=1}\frac{\sum_{d|n}A(d)(-d)^{\nu}\mu\left(\frac{n}{d}\right)}{e^{nx}-1}
\end{equation}  
\textbf{Proof.}\\
$$
\frac{d^{\nu}}{dx^{\nu}}\left(\sum^{\infty}_{n=1}\frac{\sum_{d|n}A(d)\mu\left(\frac{n}{d}\right)}{e^{nx}-1}\right)=\frac{d^{\nu}}{dx^{\nu}}f\left(e^{-x}\right)=\sum^{\infty}_{n=1}\frac{f^{(n)}(0)}{n!}\frac{d^{\nu}}{dx^{\nu}}\left(e^{-nx}\right)=
$$
$$
=\sum^{\infty}_{n=1}\frac{f^{(n)}(0)}{n!}(-n)^{\nu}e^{-nx}
$$
Using again Proposition 1 we get the result.\\
\\
\textbf{Lemma 1.}\\
\begin{equation}
\sum^{\infty}_{n=1}\frac{X(n)}{e^{nx}-1}=\sum^{\infty}_{n=1}\sum_{d|n}X(d)e^{-nx}
\end{equation}
\textbf{Proof.}\\
Set $\sum_{d|n}\frac{f^{(d)}(0)}{d!}\mu\left(\frac{n}{d}\right)=X(n)$, then from Moebius inversion theorem we have $\frac{f^{(n)}(0)}{n!}=\sum_{d|n}X(d)$. Using Proposition 2 we get the result.\\
\\
\textbf{Proposition 5.}\\
Let $\sum_{d|n}X(d)=\frac{g^{(n)}(0)}{n!}$ and $|q|<1$, then for every $f$ we have
\begin{equation}
\sum^{\infty}_{n=1}\frac{q^n}{1-q^n}\sum_{d|n}X(d)f\left(\frac{n}{d}\right)=\sum^{\infty}_{n=1}g(q^n)f(n)
\end{equation}
\textbf{Proof.}\\
Let $\sum_{d|n}X(d)=\frac{g^{(n)}(0)}{n!}$, from Lemma 1 we have
$$
\sum^{\infty}_{n=1}\frac{X(n)f(m)}{e^{nmx}-1}=\sum^{\infty}_{n=1}\frac{g^{(n)}(0)}{n!}f(m)e^{-nmx}.
$$ 
Summing with respect to $m$ we have $$\sum^{\infty}_{n=1}\frac{\sum_{d|n}X(d)f\left(\frac{n}{d}\right)}{e^{nx}-1}=\sum^{\infty}_{n=1}f(n)g\left(e^{-nx}\right)$$ and the result follows.\\
\\
\textbf{Proposition 6.}\\
Let $\sum_{d|n}X(d)=\frac{g^{(n)}(0)}{n!}$, then for every $f$ and $|q|<1$ we have
\begin{equation}
\sum^{\infty}_{n=1}\frac{q^n}{1+q^n}\sum_{d|n}X(d)f\left(\frac{n}{d}\right)=\sum^{\infty}_{n=1}\left(g(q^n)-2g(q^{2n})\right)f(n)
\end{equation}
\\
\textbf{Theorem 2.}\\
Let $|q|<1$ then
\begin{equation}
\sum^{\infty}_{n=1}\frac{q^n}{1-q^n}\sum_{d|n}f(d)\phi_H\left(\frac{n}{d}\right)=\sum^{\infty}_{n=1}f(n)H(q^n)
\end{equation}
where $\phi_H(n)=\sum_{d|n}h_d\mu\left(\frac{n}{d}\right)$ and $H(x)=\sum^{\infty}_{k=1}h_kx^k$.\\
\\
\textbf{Examples.}\\
\textbf{1)}
Set $X(n)=n^{\nu}$ in (11) then
\begin{equation}
\sum^{\infty}_{n=1}\frac{n^{\nu}}{e^{nx}-1}=\sum^{\infty}_{n=1}\sigma_{\nu}(n)e^{-nx}
\end{equation}
\\
\textbf{2)}
Also setting $\frac{f^{(n)}(0)}{n!}=n^{-\nu}$
\begin{equation}
\sum^{\infty}_{n=1}\frac{\sum_{d|n}d^{-\nu}\mu\left(\frac{n}{d}\right)}{e^{nx}-1}=\textrm{Li}_{\nu}\left(e^{-x}\right)\textrm{, }x>0
\end{equation}
or the equivalent
\begin{equation}
\sum^{\infty}_{n=1}\frac{q^n}{1-q^n}\sum_{d|n}d^{-\nu}\mu\left(\frac{n}{d}\right)=\textrm{Li}_{\nu}(q)
\end{equation}
where $\textrm{Li}_{\nu}(x)=\sum^{\infty}_{n=1}\frac{x^n}{n^{\nu}}$.\\
\\
\textbf{3)} With $h_n=\delta_n$ in Theorem 2 and $f(n)\rightarrow a(n)$, we get
\begin{equation}
\sum^{\infty}_{n=1}\frac{q^n}{1-q^n}\sum_{d|n}a(d)\mu\left(\frac{n}{d}\right)=\sum^{\infty}_{n=1}a(n)q^n
\end{equation}
Differentiating with respect to $q$ and setting $q=e^{2\pi i z}$, $Im(z)>0$, we get
\begin{equation}
\sum^{\infty}_{n=1}\frac{nX(n)}{\sin^2(\pi z n)}=4\sum^{\infty}_{n=1}a(n)n q^n\textrm{, where } X(n)=\sum_{d|n}a(d)\mu\left(\frac{n}{d}\right).
\end{equation}
The case 
\begin{equation}
X(n)=\sum^{N}_{\scriptsize
\begin{array}{cc}
	k=-M\\
	k-odd
\end{array}\normalsize
}c_kn^k,
\end{equation}  
lead us to some kind of ''Eisenstein series'' (we have set $E_{k}(q):=\sum^{\infty}_{n=1}\frac{n^{k}q^n}{1-q^n}$):
\begin{equation}
-\sum^{\infty}_{n=1}\frac{X(n)n}{\sin^2(\pi z n)}=4\sum^{N}_{\scriptsize
\begin{array}{cc}
	k=-M\\
	k-odd
\end{array}\normalsize
}c_k\sum^{\infty}_{n=1}\sigma_{k}(n)nq^n=4q\sum^{N}_{\scriptsize
\begin{array}{cc}
	k=-M\\
	k-odd
\end{array}\normalsize
}c_k\frac{d}{dq}E_{k+1}(q).
\end{equation}
Hence
\begin{equation}
-\sum^{\infty}_{n=1}\frac{X(n)n}{\sin^2(\pi z n)}=4q\frac{d}{dq}\sum^{N}_{\scriptsize
\begin{array}{cc}
	k=-M\\
	k-odd
\end{array}\normalsize
}c_kE_{k+1}(q).
\end{equation}
Set 
\begin{equation}
\Lambda(s)=M\left(\sum^{\infty}_{n=1}\frac{X(n)n}{\sinh^2(\pi t n)}\right)(s),
\end{equation}
where $M(f)(s)=\int^{\infty}_{0}f(t)t^{s-1}dt$ is the Mellin transform of the function $f(t)$. Then
\begin{equation}
\Lambda(s)=(2\pi)^{-s}\Gamma(s)\sum^{\infty}_{n=1}\sum^{N}_{\scriptsize
\begin{array}{cc}
	k=-M\\
	k-odd
\end{array}\normalsize
}c_k\sigma_k(n)n^{-s}
\end{equation}
and
\begin{equation}
\Lambda_k(s)=(2\pi)^{-s}\Gamma(s)\sum^{\infty}_{n=1}\frac{\sigma_k(n)}{n^s}=(2\pi)^{-s}\Gamma(s)\zeta(s)\zeta(s-k).
\end{equation}
Then from Hecke theorem for modular forms we have 
\begin{equation}
\Lambda(s)=\sum^{N}_{\scriptsize
\begin{array}{cc}
	k=-M\\
	k-odd
\end{array}\normalsize
}c_ki^{k+1}\Lambda_k(k+1-s)
\end{equation}
If $\textbf{M}_{k+1}$ denotes the $k+1-$th space of modular forms (that is of weight $k+1$) and
\begin{equation}
\textbf{M}=\textbf{M}_{k_1+1}\oplus \textbf{M}_{k_2+1}\ldots\oplus\textbf{M}_{k_p+1},
\end{equation} 
then we can say $f(z)\in\textbf{M}$ iff $f(z)$ can be written as a sum of $p$ different weight modular forms.\\
The same thing happens and with $f$ in which
\begin{equation}
\sum^{\infty}_{n=1}\sum^{N}_{\scriptsize
\begin{array}{cc}
	k=-M\\
	k-odd
\end{array}\normalsize
}c_kn^{k+1}\sin^{-2}(\pi z n)=4q\frac{d}{dq}f(q).
\end{equation}
Then $f$ is a sum of modular forms of different weights.\\
Hence we can state that if a function is a ''mixed'' modular form
\begin{equation}
f(z)=\sum^{\infty}_{n=0}a_f(n)q^n\textrm{, }q=e^{2\pi i z}\textrm{, }Im(z)>0
\end{equation}
with 
\begin{equation}
a_f(n)=\sum_{d|n}X(d)\textrm{, where  }X(n)=c_1X_1(n)+c_2X_2(n)+\ldots+c_pX_p(n)
\end{equation} 
and
\begin{equation}
\Lambda_{\nu}(s)=(2\pi)^{-s}\Gamma(s)\sum^{\infty}_{n=1}\frac{\sum_{d|n}X_{\nu}(d)}{n^s}
\end{equation}
is such that 
\begin{equation}
\Lambda_{\nu}(s)=i^{k_{\nu}}\Lambda_{\nu}\left(k_{\nu}-s\right)\textrm{, }\nu=1,2,\ldots,p
\end{equation}
then
\begin{equation}
\Lambda_f(s)=\sum^{p}_{\nu=1}c_{\nu}i^{k_{\nu}}\Lambda_{\nu}(k_{\nu}-s),
\end{equation}
where
\begin{equation}
\Lambda_f(s)=(2\pi)^{-s}\Gamma(s)\sum^{\infty}_{n=1}\frac{a_f(n)}{n^s}
\end{equation}
and
\begin{equation}
M_f(s)=\frac{\pi}{2(s-1)}\int^{\infty}_{0}\left(\sum^{\infty}_{n=1}\frac{X(n)n}{\sinh^2(\pi n t)}\right)t^{s-1}dt=\Lambda_f(s-1). 
\end{equation}
\\

Moreover holds the next Hecke-type theorem for derivatives of modular forms:\\ 
\\
\textbf{Theorem 3.}\\
If $f(z)$ is modular form of even weight $k$ and we set
\begin{equation}
g(z):=\sum^{\infty}_{n=1}\frac{X(n)n}{\sin^2(\pi n z)}=4q\frac{d}{dq}f(z)\textrm{, }a_f(n)=\sum_{d|n}X(d),
\end{equation}
then we have
\begin{equation}
\Lambda^{*}(s)=\int^{\infty}_{0}g(it)t^{s-1}dt=4\frac{\Gamma(s)}{(2\pi)^s}\sum^{\infty}_{n=1}\frac{a_f(n)n}{n^s}=\frac{2}{\pi}(s-1)\Lambda_f(s-1),
\end{equation}
where $\Lambda_f(s)=(2\pi)^{-s}\Gamma(s)\sum^{\infty}_{n=1}\frac{a_f(n)}{n^s}$ and $\Lambda^{*}(s)$ satisfies the functional equation
\begin{equation}
\frac{\Lambda^{*}(s)}{s-1}=i^{k}\frac{\Lambda^{*}(k+2-s)}{k+1-s}.
\end{equation}
\\
\textbf{Remarks.}
i) Here the weight is $k$. 
ii) The proof of the theorem is based on Proposition 11 and Hecke's theorem.\\

Assume that 
\begin{equation}
f(z)=\sum^{\infty}_{n=0}a_f(n)q^n\textrm{, }q=e^{2\pi i z}=e(z)\textrm{, }Im(z)>0.
\end{equation}
Set also
\begin{equation}
q\frac{d}{dq}f(z)=\sum^{\infty}_{n=1}c_f(n)q^n
\end{equation}
then 
\begin{equation}
c_f(n)=a_f(n)n
\end{equation}
and hence
\begin{equation}
c_f(nm)\equiv0 \textrm{mod} nm\textrm{, for all }n,m\in\{1,2,\ldots\}
\end{equation}
\begin{equation}
c_f(n)c_f(m)\equiv 0\textrm{mod} nm\textrm{, for all }n,m\in\{1,2,\ldots\}
\end{equation}
From the above we get the following\\ 
\\
\textbf{Lemma 1.}\\
If $a_f(n)$ are integers, then
\begin{equation}
c_f(nm)\equiv c_f(n)c_f(m)\textrm{mod}nm
\end{equation}
\\
\textbf{Lemma 2.}\\
Suppose that $f(q)$ is differentiatable function of $q$, $|q|<1$, then $f(e(z))$, $z=x+iy$, $x,y\in \textbf{R}$ is harmonic, in the sence $h(x,y):=f(e(x+iy))$ satisfies the equation
\begin{equation}
\partial^2_xh(x,y)+\partial^2_yh(x,y)=0
\end{equation}
\\
\textbf{Lemma 3.}\\
There exists function $K(x,y)\in C^{\infty}_0(\textbf{R}^2)$ such that
\begin{equation}
f(e(z))=\int^{+\infty}_{-\infty}\int^{+\infty}_{-\infty}\frac{yK(x',y')}{\left(y^2+(x-x')^2+(y-y')^2\right)^{3/2}}dx'dy'
\end{equation}

\section{Results in the theory of theta functions}

Set $(n,m)=gcd(n,m)$ to be the greatest common divisor of $n,m$. Then one can easily see, using arguments of [1] chapter 2, that 
\begin{equation}
\sum^{n}_{m=1}f((n,m))=\sum_{d|n}f(d)\phi\left(\frac{n}{d}\right).
\end{equation} 
From the relation $\sum_{d|n}\phi(d)=n$ and relation (11) we get
\begin{equation}
\sum^{\infty}_{n=1}\frac{\phi(n)}{e^{nx}+1}=\frac{1}{2}\frac{\cosh(x)}{\sinh^2(x)}
\end{equation}
From Proposition 5 we have for general arithmetic function $F$: 
$$\sum^{\infty}_{n=1}\frac{\sum_{d|n}\phi(d)F\left(\frac{n}{d}\right)}{e^{nx}+1}=\frac{1}{2}\sum^{\infty}_{n=1}F(n)\frac{\cosh(nx)}{\sinh^2(nx)}$$
or else
\begin{equation}
2\sum^{\infty}_{n=1}\frac{\sum^{n}_{m=1}F((n,m))}{e^{nx}+1}=\sum^{\infty}_{n=1}F(n)\frac{\cosh(nx)}{\sinh^2(nx)}
\end{equation}
integrating the above relation we get:\\
\\
\textbf{Theorem 4.}
\begin{equation}
\prod^{\infty}_{n=1}(1+q^n)^{\frac{1}{n}\sum^{n}_{m=1}F((n,m))}=\exp\left(\sum^{\infty}_{n=1}\frac{F(n)q^n}{n(1-q^{2n})}\right)
\end{equation}

Set now $f(-q):=\prod^{\infty}_{n=1}(1-q^n)$, then we have the next version of Jacobi triple product identity (see [4]):
\begin{equation}
\sum^{\infty}_{n=1}\frac{\cosh(tn)}{n\sinh(\pi a n)}=\log(f(-e^{-2\pi a}))-\log\left(\theta_4(it/2,e^{-a\pi})\right)\textrm{, }|t|<a\pi
\end{equation}
Setting $e^{-a\pi}=q$, $a>0$ and $F(x)=\cos(2tx)$ in (50) and using (51) we get\\
\\
\textbf{Proposition 7.}\\
If $|q|<1$, then 
\begin{equation}
\prod^{\infty}_{n=1}(1+q^n)^{\frac{1}{n}\sum^{n}_{m=1}\cos(2t(n,m))}=\left(\frac{f(-q^2)}{\theta_4(t,q)}\right)^{1/2},
\end{equation}
where 
$$
\theta_4(z,q)=\sum^{\infty}_{n=-\infty}(-1)^nq^{n^2}e^{2niz}\textrm{, }|q|<1
$$ 
is Jacobi's 4th theta function (see [2]).\\  
\\

Setting $t=\frac{\pi}{2}$ in (52) we get
\begin{equation}
q^{1/2}\prod^{\infty}_{n=1}(1+q^n)^{\frac{12}{n}\sum^{n}_{m=1}(-1)^{(n,m)}}=q^{1/2}\left(\frac{f(-q^2)}{\theta_4\left(\frac{\pi}{2},q\right)}\right)^6=q^{1/2}\left(\frac{f(-q^2)}{\theta_3\left(q\right)}\right)^6,
\end{equation}
where $\theta_3(q)$ is the ''null'' theta function $\theta_4\left(\frac{\pi}{2},q\right)=\sum^{\infty}_{n=-\infty}q^{n^2}$.\\Hence we can state the next\\ 
\\
\textbf{Proposition 8.}\\
If $q=e^{-\pi\sqrt{r}}$, $r>0$, then
\begin{equation}
q^{1/2}\prod^{\infty}_{n=1}(1+q^n)^{\frac{12}{n}\sum^{n}_{m=1}(-1)^{(n,m)}}=4kk'.
\end{equation} 
\\
\textbf{Proof.}\\
We use (see [6] pg.488):
\begin{equation}
f(-q^2)^6=\prod^{\infty}_{n=1}\left(1-q^{2n}\right)^6=\frac{2kk'K(k)^3}{\pi^3q^{1/2}}
\end{equation}
and (see [2] pg.107 and related theory): 
\begin{equation}
\theta_3(q)^2=\left(\sum^{\infty}_{n=-\infty}q^{n^2}\right)^2=\frac{2K}{\pi}
\end{equation}
and relation (53). The functions $k=k_r$, $K=K(k)$ are the elliptic singular modulus and elliptic integral of the first kind at singular values respectively. The function $k'=\sqrt{1-k^2}$ is called complementary modulus (see [2] pg.11).\\

Also if $q=e^{-\pi\sqrt{r}}$, $r>0$ we set
\begin{equation}
\psi(\nu,z,q)=\sum^{\infty}_{n=1}\frac{n^{\nu}z^n}{n(1-q^n)}.
\end{equation} 
Then using the summable version of Jacobi triple product identity (relation (51)), we get after differentiating with respect to $t$:
\begin{equation}
\psi(2\nu,q,q^2)=\sum^{\infty}_{n=1}\frac{n^{2\nu}q^n}{n(1-q^{2n})}=-\frac{1}{2(-4)^{\nu}}\left[\frac{\partial ^{2\nu}}{\partial t^{2\nu}}\log\left(\theta_4(t,q)\right)\right]_{t=0}.
\end{equation}
Hence from (50) we have\\
\\
\textbf{Theorem 5.}\\
If $q=e^{-\pi\sqrt{r}}=e^{i\pi\tau}$, $\tau=\sqrt{-r}$, $r>0$, then
$$
\log\left(\prod^{\infty}_{n=1}\left(1+q^n\right)^{\frac{1}{n}\sum^{n}_{m=1}(n,m)^{2\nu}}\right)=\frac{(-1)^{\nu+1}}{2^{2\nu+1}}\left[\frac{\partial^{2\nu}}{\partial t^{2\nu}}\log\left(\theta_4(t,q)\right)\right]_{t=0}=
$$
\begin{equation}
=\frac{B_{2\nu}}{4\nu}\left(E_{2\nu}(q^2)-E_{2\nu}(q)\right)
\end{equation}
is a modular form of weight $2\nu$,
where
\begin{equation}
E_{\nu}(q)=1-\frac{2\nu}{B_{\nu}}\sum^{\infty}_{n=1}\frac{n^{\nu-1}q^n}{1-q^n}
\end{equation}
are the known Eisenstein series (see [18]).\\
\\
\textbf{Theorem 6.}\\
For $\nu=1,2,\ldots$, we have
\begin{equation} 
\sum^{\infty}_{n=1}\frac{n^{2\nu-1}}{\sinh(\pi n\sqrt{r})}=\frac{B_{2\nu}}{2\nu}\left(E_{2\nu}\left(q^2\right)-E_{2\nu}(q)\right)\textrm{, }q=e^{-\pi\sqrt{r}}\textrm{, }r>0
\end{equation}
\\
\textbf{Proof.}\\
Easy. From (60) we immediately get the result.\\
\\
\textbf{Theorem 7.}\\
Let $F(n)=a_2n^2+a_4n^4+\ldots+a_{2k}n^{2k}+\ldots$ be even function with $F(0)=0$, then
\begin{equation}
\prod^{\infty}_{n=1}\left(1+q^n\right)^{\frac{1}{n}\sum^{n}_{m=1}F\left((n,m)\right)}=\exp\left(\sum^{\infty}_{\nu=1}\frac{a_{2\nu}B_{2\nu}}{4\nu}\left(E_{2\nu}(q^2)-E_{2\nu}(q)\right)\right).
\end{equation}
\\

Another interesting result which follows from relation 
(59) is
\begin{equation}
\log\left(\frac{\theta_4(t,q)^{1/2}\theta_4(-t,q)^{1/2}}{\theta_4(q)}\right)=\log\left(\frac{\theta_4(t,q)}{\theta_4(q)}\right)=-2\sum^{\infty}_{n,m=1}\frac{(-4)^n}{(2n)!}B(n,m)t^{2n}q^m
\end{equation}
where $B(\nu,n)=-\frac{1}{n}\sum_{d|n}(-1)^{n/d}\sum^{d}_{k=1}(d,k)^{2\nu}$.\\ This can be done by writing 
$$
\sum^{\infty}_{n=1}\left({\frac{1}{n}\sum^{n}_{m=1}(n,m)^{2\nu}}\right)\log\left(1+q^n\right)=\frac{(-1)^{\nu+1}}{2^{2\nu+1}}\left[\frac{\partial^{2\nu}}{\partial t^{2\nu}}\log\left(\theta_4(t,q)\right)\right]_{t=0}
$$
Expanding the logarithm on the left into power series we get
$$
\frac{(-1)^{\nu+1}}{2^{2\nu+1}}\left(\frac{\partial^{2\nu}}{\partial t^{2\nu}}\log(\theta_4(t,q))\right)_{t=0}=-\sum^{\infty}_{n=1}q^{n}\frac{1}{n}\sum_{d|n}(-1)^{n/d}\sum^{d}_{m=1}(d,m)^{2\nu}
$$
From this we arrive to (66) and we have the next\\
\\
\textbf{Proposition.}\\
If $|q|<1$, then
\begin{equation}
\log\left(\frac{\theta_4(t,q)}{\theta_4(q)}\right)=-2\sum^{\infty}_{n=1}P(n,t)q^n ,
\end{equation}
with $P(n,t)=\frac{1}{n}\sum_{d|n}(-1)^{n/d}\sum^{d}_{m=1}\left(1-\cos\left(2(d,m)t\right)\right)$.\\
\\
\textbf{Notes.}\\
Relation (67) can also be found if we use directly the Jacobi's formula
\begin{equation}
\frac{\partial_z\theta_4(z,q)}{\theta_4(z,q)}=4\sum^{\infty}_{n=1}\frac{q^n}{1-q^{2n}}\sin(2nz).
\end{equation} 
This can be done integrating (68) from $z=0$ to $z=t$. Then we arive to (67) expanding $\frac{q^n}{1-q^{2n}}=\frac{x}{1-x^2}$ into Taylor series of $x=q^n$ and making the double infinite series into single using divisor sumation.\\

The above formula (67) for $t=\frac{\pi}{2}$ gives\\
\\
\textbf{Proposition 9.}\\
Let $|q|<1$, then
\begin{equation}
\log\left(\frac{\theta_3(q)}{\theta_4(q)}\right)=-\frac{1}{2}\log(k')=\sum^{\infty}_{n=1}P_1(n)q^n
\end{equation}
where $P_1(n)=\frac{1}{n}\sum_{d|n}(-1)^{n/d}\sum^{d}_{k=1}\left(1-(-1)^{(d,k)}\right)$.\\ 
\\
\textbf{Definition 1.}\\
We define $(Ta)(n)$ to be such that
\begin{equation}
\exp\left(-\sum^{\infty}_{n=1}a_nx^n\right)=\sum^{\infty}_{n=0}(Ta)(n)x^n,
\end{equation}
\\
\textbf{Theorem 9.}\\
The above transformation $Ta(n)=b_n$ is given from relation
\begin{equation}
nb_n=-\sum^{n}_{k=1}a_kkb_{n-k}\textrm{, }n=1,2,\ldots ,
\end{equation}
where $b_0=1$. Hence given the general form of $a_n$, we can construct all $Ta(n)$.\\
\\
\textbf{Notes.} In general hold the following relations\\ 
\textbf{1.} Set $c_n=a_n+b_n$, then 
\begin{equation}
(Tc)(n)=\sum^{n}_{k=0}(Ta)(n-k)(Tb)(k)
\end{equation}
\textbf{2.} If $c_n=a_n+\frac{1}{n}$, then
\begin{equation}
(Tc)(n)=(Ta)(n)-(Ta)(n-1)
\end{equation}
\textbf{3.} Consider the function
\begin{equation}
\exp\left(\frac{x}{x-1}\right)(1-x)=\sum^{\infty}_{n=0}a(n)x^n,
\end{equation}
then 
\begin{equation}
\sum^{\infty}_{n=0}a(n)=0.
\end{equation}
This happens because
\begin{equation}
\exp\left(-\sum^{\infty}_{n=1}x^n\right)=\exp\left(\frac{x}{x-1}\right)=\sum^{\infty}_{n=0}h(n)x^n
\end{equation}
and $a(n)=(T1)(n)-(T1)(n-1)=h(n)-h(n-1)$, with $\lim_{n\rightarrow\infty}h(n)=0$. Infact it is 
\begin{equation}
\exp\left(-\sum^{\infty}_{n=1}x^n\right)=\exp\left(\frac{x}{x-1}\right)=\sum^{\infty}_{n=0}h(n)x^n
\end{equation}
and
\begin{equation}
\sum^{N}_{n=1}h(n)=(T1)(N-1),
\end{equation}
with $\lim_{N\rightarrow\infty}(T1)(N)=0$.\\

Continuing we set the numbers $s(n)$ such that
\begin{equation}
\theta_4(q)=\exp\left(\sum^{\infty}_{n=1}s(n)q^n\right),
\end{equation}
then (our thoughts motivated from [16]):
\begin{equation}
n\sum_{d|n}s(d)\mu\left(\frac{n}{d}\right)=1-(-1)^n
\end{equation}
Hence using Moebius inversion theorem
\begin{equation}
s(n)=\sum_{d|n}\frac{1-(-1)^d}{d}.
\end{equation}
\\
From (67), the property of $s(n)$ and Definition 1, we get\\
\\
\textbf{Theorem 10.}\\
For $|q|<1$ and $t$ real, we have 
\begin{equation}
\theta_4(t,q)=\exp\left(-\sum^{\infty}_{n=1}\left(2P(n,t)+\sum_{d|n}\frac{(-1)^d-1}{d}\right)q^n\right).
\end{equation}
and
\begin{equation}
\theta_4(t,q)=\sum^{\infty}_{n=0}T\left(2P(n,t)+\sum_{d|n}\frac{(-1)^d-1}{d}\right)q^n.
\end{equation}
\\

If we consider the theta function $
\theta_3(q)=\sum^{\infty}_{n=-\infty}q^{n^2}\textrm{, }|q|<1
$, then writing
\begin{equation}
\theta_3(q)=\exp\left(\sum^{\infty}_{n=1}t(n)q^n\right),
\end{equation}
we get
\begin{equation}
n\sum_{d|n}t(d)\mu\left(\frac{n}{d}\right)=-\psi_4(n),
\end{equation}
where $\psi_4$ is a 4-periodic arithmetic function with
\begin{equation}
\psi_4(n)=\left\{
\begin{array}{cc}
   0\textrm{ if }n\equiv 0(mod4)\\
-2\textrm{ if }n\equiv 1(mod4)\\
	 8\textrm{ if }n\equiv 2(mod4)\\
-2\textrm{ if }n\equiv 3(mod4)
\end{array}\right\}
.\end{equation}
Hence as above\\
\\
\textbf{Theorem 11.}\\
The next expansion is valid
\begin{equation}
\theta_3(q)=\exp\left(-\sum^{\infty}_{n=1}\sum_{d|n}\frac{\psi_4(d)}{d}q^n\right)\textrm{, }|q|<1
\end{equation} 
and
\begin{equation}
r_{\nu}(n)=T\left(\nu\sum_{d|n}\frac{\psi_4(d)}{d}\right)(n),
\end{equation}
where $r_{\nu}(n)$ is the number of the representations of the positive integer $n$ into $\nu$ squares.\\ 

We can use the Chebyshev polynomials $C_n(x)$ to recover the theta functions. From the property 
\begin{equation}
C_n(\cos(x))=\cos(nx),
\end{equation}
we get\\
\\
\textbf{Proposition 10.}\\
If $|q|<1$, then
\begin{equation}
\log\left(\frac{\theta_4(\arccos(t),q)}{\theta_4(q)}\right)=-2\sum^{\infty}_{n=1}P^{*}(n,t)q^n
\end{equation}
where $P^{*}(n,t)=\frac{1}{n}\sum_{d|n}(-1)^{n/d}\sum^{d}_{m=1}(1-C_{2(d,m)}(t))$ and $C_n(x)$ is the $n$-th order Chebyshev orthogonal polynomial, hence the $P^{*}(n,t)$ is a $2n$-degree even polynomial in $t$.\\
\\
\textbf{Proposition.}\\
Let $f$ have Fourier series expansion  
\begin{equation}
f(\phi)=\sum^{\infty}_{n=1}a_{n}\cos(2n\phi),
\end{equation}
then
\begin{equation}
\int^{\pi}_{0}f(\phi)\log\left(\frac{\theta_4(\phi,q)}{\theta_4(q)}\right)d\phi=-\pi\sum^{\infty}_{n=1}\frac{a_{n}q^n}{n(1-q^{2n})}
\end{equation}
\textbf{Proof.}\\
From Proposition 10 and using the orthogonality of Chebyshev polynomials we get 
\begin{equation}
\int^{1}_{-1}\frac{f(t)}{\sqrt{1-t^2}}\log\left(\frac{\theta_4(\arccos(t),q)}{\theta_4(q)}\right)dt=-\pi\sum^{\infty}_{k=1}\sum^{\infty}_{n=1}P_1(k,n)a_{2k}q^{n/2} ,
\end{equation}
where $P_1(k,n)=1/k$ if $k|n$ while $2k\neq0(\textrm{mod} n)$ and otherwise 0.\\ 
The function $P_1$ is a characteristic periodic arithmetical function and can be evaluated using arithmetic progressions. After summing with respect to $n$ and making the change of variable $t\rightarrow \cos(\phi)$ in the above integral we get the result.\\
\\
\textbf{Notes.}\\
\textbf{1)} The author feels to ask if the above problem may be treated with countor integration and residuals calculus, since the following very interesting formula rises from the work of Chouika (see [20]): 
\begin{equation}
\theta_4(u,q)=\theta_4(q)\prod^{\infty}_{n=0}\left(1-\left(\frac{\sin(\pi u)}{\sin((n+1/2)\pi\tau)}\right)^2\right)\textrm{, }q=e^{\pi i \tau}\textrm{, }Im(\tau)>0 
\end{equation} 
\textbf{2)} For to get integral (92), we pass through   
\begin{equation}
\int^{1}_{-1}\frac{f(t)}{\sqrt{1-t^2}}\log\left(\frac{\theta_4(\arccos(t),q)}{\theta_4(q)}\right)dt=-\pi\sum^{\infty}_{n=1}\frac{a_{2n}q^n}{n(1-q^{2n})},
\end{equation}
when $f(x)=\sum^{\infty}_{n=1}a_{2n}C_{2n}(x)$.\\
\\
\textbf{Examples.}\\
\textbf{i)} For $f(x)=C_2(x)=-1+2x^2$ we get 
\begin{equation}
\int^{1}_{-1}\frac{1-2t^2}{\sqrt{1-t^2}}\log\left(\frac{\theta_4(\arccos(t),q)}{\theta_4(q)}\right)dt=\frac{\pi q}{1-q^2}
\end{equation}
\textbf{ii)} 
If $a_n=n^{-3}$ then $f(x)=8^{-1}\sum^{\infty}_{n=1}\frac{\cos(2n\arccos(x))}{n^3}$ and hence if $q=e^{-x}$, $x>0$ we get
\begin{equation}
\int^{1}_{-1}\frac{f(t)}{\sqrt{1-t^2}}\log\left(\frac{\theta_4(\arccos(t),q)}{\theta_4(q)}\right)dt=-\frac{\pi}{16}\sum^{\infty}_{n=1}\frac{1}{n^4\sinh(nx)}
\end{equation}
\textbf{iii)} If $a_n=\frac{1}{n}$ then $f(x)=-\frac{1}{2}\log(2\sin(x))$ and hence
\begin{equation}
\int^{\pi}_{0}\log(2\sin(\phi))\log\left(\frac{\theta_4(\phi,q)}{\theta_4(q)}\right)d\phi=\pi\sum^{\infty}_{n=1}\frac{q^n}{n^2(1-q^{2n})}
\end{equation}
\textbf{iv)} If $a_n=q^{n^2}$, then $f(x)=\sum^{\infty}_{n=1}q^{4n^2}\cos(2nx)=\frac{\theta_3(x,q^4)-1}{2}$
\begin{equation} 
\int^{\pi}_{0}\frac{\theta_3(\phi,q^4)-1}{2}\log\left(\frac{\theta_4(\phi,q)}{\theta_4(q)}\right)d\phi=-\pi\sum^{\infty}_{n=1}\frac{q^{4n^2+n}}{n(1-q^{2n})}
\end{equation}
\textbf{v)} If $a_n=nq^{n^2/2+n/2}$, then 
\begin{equation} 
\int^{\pi}_{0}\psi(\phi,q^{1/2})\log\left(\frac{\theta_4(\phi,q^{1/2})}{\theta_4(q^{1/2})}\right)d\phi=-2\pi\sum^{\infty}_{n=1}\frac{q^{n^2+n}}{1-q^{n}}
\end{equation}
where $\psi(\phi,q)=2\sum^{\infty}_{n=1}nq^{2n^2+n}\cos(2n\phi)$, where $|q|<1$.\\

In the same way as in (20),(21) we obtain 
\begin{equation}
2\sum^{\infty}_{n=1}\frac{\sum^{n}_{m=1}F((n,m))}{e^{nx}-1}=\sum^{\infty}_{n=1}\frac{F(n)}{\cosh(nx)-1}
\end{equation}
Taking in both sides of (100) the Mellin transform (the Mellin transform of the function $g(x)$ is $M(g(x))(s)=M(g)(s)=\int^{\infty}_{0}g(t)t^{s-1}dt$), we get 
$$
2\Gamma(s)\zeta(s-1)\sum^{\infty}_{n=1}\frac{F(n)}{n^s}=\sum^{\infty}_{n=1}F(n)\int^{\infty}_{0}\frac{x^{s-1}}{\cosh(nx)-1}dx.
$$
Hence if $F(n)=X(n)$ is arithmetic function, in view of [1] we get\\
\\
\textbf{Proposition 11.}
\begin{equation}
M\left(\sum^{\infty}_{n=1}\frac{X(n)}{\sinh^2(nx)}\right)(s)=4\cdot2^{-s}\Gamma(s)\zeta(s-1)L(X,s)
\end{equation}
where $L(X,s)=\sum^{\infty}_{n=1}\frac{X(n)}{n^s}$. 
In case $X(n)=X_T(n)$ is $T$ periodic, then
\begin{equation}
M\left(\sum^{\infty}_{n=1}\frac{X_T(n)}{\sinh^2(nx)}\right)(s)=4\cdot 2^{-s}\Gamma(s)\zeta(s-1)T^{-s}\sum^{T}_{m=1}X_T(m)\zeta\left(s,\frac{m}{T}\right)
\end{equation}
The function $\zeta(z,\nu)=\sum^{\infty}_{n=0}\frac{1}{(z+n)^{\nu}}$ is the Hurwitz zeta function and is a generalization of the Riemann's $\zeta(\nu)$ function.\\
\\
\textbf{Proposition 12.}(see [4])\\
If 
\begin{equation}
X_1(n)= \left\{\begin{array}{cc}
 1, \mbox{  }  n\equiv p-a(modp) \\
-1, \mbox{  }  n\equiv p-b (modp) \\
 1, \mbox{  }  n\equiv a(modp) \\
-1, \mbox{  }  n\equiv b(modp) \\
 0, \mbox{  }  else 
\end{array}\right\},
\end{equation}  
then 
\begin{equation}
\prod^{\infty}_{n=1}(1-q^n)^{X_1(n)}=\frac{\vartheta_4\left((p-2a)ix/4,e^{-px/2}\right)}{\vartheta_4\left((p-2b)ix/4,e^{-px/2}\right)} 
\end{equation}
\begin{equation}
\sum^{\infty}_{n=1}\frac{X_1(n)n^2}{\sinh^2(nx)}=-\frac{d^2}{dx^2}\log\left(\frac{\vartheta_4\left((p-2a)ix/2,e^{-px}\right)}{\vartheta_4\left((p-2b)ix/2,e^{-px}\right)}\right)
\end{equation}
and
\begin{equation}
\int^{\infty}_{0}\log\left(\frac{\vartheta_4\left((p-2a)ix/4,e^{-px/2}\right)}{\vartheta_4\left((p-2b)ix/4,e^{-px/2}\right)}\right)x^{s-1}dx=-\frac{\Gamma(s)\zeta(s+1)}{p^s}\sum^{p}_{m=1}X_1(m)\zeta(s,\frac{m}{p})
\end{equation}
\\
\textbf{Examples.}\\
\textbf{i)} If $q=e^{-2x}$, then
\begin{equation}
\sum^{\infty}_{n=1}\left
(\frac{n}{8}\right)\frac{n^2}{\sinh^2(nx)}=-\frac{d^2}{dx^2}\log\left(\frac{\sum^{+\infty}_{n=-\infty}(-1)^nq^{4n^2+3n}}{\sum^{+\infty}_{n=-\infty}(-1)^nq^{4n^2+n}}\right)
\end{equation}
\textbf{ii)} If $q=e^{-2x}$, then
\begin{equation}
\sum^{\infty}_{n=1}\left
(\frac{n}{5}\right)\frac{n^2}{\sinh^2(nx)}=-\frac{d^2}{dx^2}\log\left(\frac{\sum^{+\infty}_{n=-\infty}(-1)^nq^{5n^2/2+3n/2}}{\sum^{+\infty}_{n=-\infty}(-1)^nq^{5n^2/2-n/2}}\right)
\end{equation}
\textbf{iii)} If $|q|<1$, then
\begin{equation}
\prod^{\infty}_{n=1}\left(1-q^n\right)^{\left(\frac{n}{5}\right)}=\frac{1}{1+}\frac{q}{1+}\frac{q^2}{1+}+\ldots
\end{equation}
\\
\textbf{Proposition 13.}\\
In general if $q=e^{-2x}$ and $X(n)$ is arithmetic function then
\begin{equation}
\sum^{\infty}_{n=1}X(n)\frac{n^2}{\sinh^2(nx)}=-\frac{d^2}{dx^2}\log\left(\prod^{\infty}_{n=1}\left(1-e^{-2nx}\right)^{X(n)}\right)
\end{equation}
\\
\textbf{Examples.}\\
\\
\textbf{i)}
\begin{equation}
\sum^{\infty}_{n=1}\frac{\mu(n)n}{\sinh^2(nx)}=4e^{-2x}
\end{equation}
\\
\textbf{ii)}
\begin{equation}
\sum^{\infty}_{n=1}\frac{n^2}{\sinh^2(nx)}=\frac{d^2}{dx^2}\log\left(\eta(2x)\right),
\end{equation}
where $\eta(x)=\prod^{\infty}_{n=1}\left(1-e^{-nx}\right)$ is the Ramanujan-Dedekind eta function.\\
\\
\textbf{iii)} If $G=g_1^{a_1}g_2^{a_2}\ldots g_s^{a_s}$, with $a_1,a_2,\ldots,a_s$ non negative integers is the prime decomposition of a perfect square, then (see [4],[16])
$$
\sum^{\infty}_{n=1}\left(\frac{n}{G}\right)\frac{n^2}{\sinh^2(nx)}=
$$
\begin{equation}
=-\frac{d^2}{dx^2}\log\left(\eta(2x)\prod^{\lambda}_{i=1}\eta(2g_ix)^{-1}\prod_{i<j}\eta(2g_ig_jx)^1\prod_{i<j<k}\eta(2g_ig_jg_kx)^{-1}\ldots\right)
\end{equation}

For example
\begin{equation}
\sum^{\infty}_{n=1}\left(\frac{n}{25}\right)\frac{n^2}{\sinh^2(nx)}=-\frac{d^2}{dx^2}\log\left(\frac{\eta(2x)}{\eta(10x)}\right).
\end{equation}
\\

From [16] and (110) we have the next:\\ 
\\
\textbf{Theorem 13.}(Conjecture)\\
If $G=2^a g_1^{a_1}g_2^{a_2}\ldots g_s^{a_s}$, with $a_1,a_2,\ldots,a_s$ non negative integers, $a\neq1$ and $g_1<g_2<\ldots<g_s$ are primes congruent to $1(\textrm{mod}4)$, then
$$
\sum^{\infty}_{n=1}\left(\frac{n}{G}\right)\frac{n^2}{\sinh^2(nx)}=-\frac{d^2}{dx^2}\log\left(\prod^{\left[\frac{G-1}{2}\right]}_{j=1}\vartheta\left(\frac{G}{2},\frac{G}{2}-j;e^{-2x}\right)^{\left(\frac{j}{G}\right)}\right).
$$
Also 
\begin{equation}
\prod^{\infty}_{n=1}(1-q^n)^{\left(\frac{n}{G}\right)}=\prod^{\left[\frac{G-1}{2}\right]}_{j=1}\vartheta\left(\frac{G}{2},\frac{G}{2}-j,q\right)^{\left(\frac{j}{G}\right)}
\end{equation}
and
\begin{equation}
\sum^{\infty}_{n=1}\left(\frac{n}{G}\right)\frac{n}{\sinh(nx)}=2\left\{q\frac{d}{dq}\log\left(\prod^{\left[\frac{G-1}{2}\right]}_{j=1}\vartheta\left(\frac{G}{2},\frac{G}{2}-j;q\right)^{\left(\frac{j}{G}\right)}\right)\right\}^{q=e^{-2x}}_{q=e^{-x}},
\end{equation}
where 
\begin{equation}
\vartheta(k,l;q)=\sum^{\infty}_{n=-\infty}(-1)^nq^{kn^2+ln}.
\end{equation}
\\
\textbf{Examples.}
\begin{equation}
\sum^{\infty}_{n=1}\left(\frac{n}{13}\right)\frac{n^2}{\sinh^2(nx)}=-\frac{d^2}{dx^2}\log\left(\prod^{4}_{j=1}\vartheta\left(\frac{13}{2},\frac{13}{2}-j;e^{-2x}\right)^{\left(\frac{j}{13}\right)}\right).
\end{equation}
and
$$
2^{-1}\sum^{\infty}_{n=1}\left(\frac{n}{5}\right)\frac{n}{\sinh(nx)}=\left\{q\frac{d}{dq}\log\left(q^{-1/5}R(q)\right)\right\}^{q=e^{-2x}}_{q=e^{-x}}=
$$
$$
=e^{-2x}\frac{R'\left(e^{-2x}\right)}{R\left(e^{-2x}\right)}-e^{-x}\frac{R'\left(e^{-x}\right)}{R\left(e^{-x}\right)}.
$$
But it is known that
\begin{equation}
\frac{R'(q)}{R(q)}=5^{-1}\frac{f\left(-q\right)^5}{qf\left(-q^5\right)}
\end{equation}
The above equality (120) has given by Ramanujan and proved later by Andrews [21] althought in the more detailed version
\begin{equation}
R(q)=\frac{\sqrt{5}-1}{2} \exp\left(\frac{1}{5}\int^{q}_{1}\frac{f^5(-t)}{tf(-t^5)}dt\right).
\end{equation}
Hence
\begin{equation}
\frac{5}{2}\sum^{\infty}_{n=1}\left(\frac{n}{5}\right)\frac{n}{\sinh(nx)}=\frac{f\left(-q^2\right)^5}{f\left(-q^{10}\right)}-\frac{f\left(-q\right)^5}{f\left(-q^5\right)}
\end{equation}

Let now
\begin{equation} 
\chi_0(n)=\left\{
\begin{array}{cc}

	-2 \textrm{ if } n\equiv

	1(mod4)
	\\
	 
	 3 \textrm{ if } n\equiv 
	 
	 2(mod4)
\\	 
	 
	-2 \textrm{ if } n\equiv 
	
	3(mod4)
\\	
	
	 1 \textrm{ if } n\equiv 
	 
	 0(mod4)	 
\end{array}
\right\} 
\end{equation}
then if $q=e^{-2x}$ we get (see [15])
\begin{equation}
\sum^{\infty}_{n=1}\chi_0(n)\frac{n^2}{\sinh^2(n x)}=-\frac{d^2}{dx^2}\log\left(\sum^{\infty}_{n=-\infty}q^{n^2}\right)=-\frac{d^2}{dx^2}\log\left(\theta_3\left(e^{-2x}\right)\right)
\end{equation} 
and
\begin{equation}
\sum^{\infty}_{n=1}\frac{(-1)^nn^2}{\sinh^2(n x)}=-\frac{d^2}{dx^2}\log\left(\sum^{\infty}_{n=-\infty}q^{n(n+1)/2}\right)
\end{equation} 
The general correspondence of hyperbolic sine function series and theta functions is the following:\\
\\
\textbf{Theorem 14.}\\
If $q=e^{-2x}$, $x>0$ and
\begin{equation} 
\chi_{k,h}(n):=\left\{
\begin{array}{cc}
  
	1 \textrm{ if } n\equiv 
	 
	0,k+h,k-h(mod2k)
\\	
	
	0 \textrm{ otherwise }	 

\end{array}
\right\}, 
\end{equation}
then
\begin{equation}
\sum^{\infty}_{n=1}\frac{\chi_{k,h}(n)n^2}{\sinh^2(n x)}=-\frac{d^2}{dx^2}\log\left(\sum^{\infty}_{n=-\infty}(-1)^nq^{kn^2+hn}\right),
\end{equation} 
when $k>h$, $k\in\bf N\rm$, $h\in\bf Z\rm$ and
\begin{equation}
\log\left(\sum^{\infty}_{n=-\infty}(-1)^nq^{kn^2+hn}\right)=-\sum^{\infty}_{n=1}f_{k,h}(n)q^n,
\end{equation}
where 
\begin{equation}
f_{k,h}(n):=\frac{1}{n}\sum_{d|n}\chi_{k,h}(d)d.
\end{equation}
\\
\textbf{Proof.}\\
The proof is similar to the proof of the next theorem.\\
\\ 
\textbf{Theorem 15.}\\
If $|q|<1$, $k>|h|>0$, $k\in\textbf{N}$, $h\in\textbf{Z}$, then
\begin{equation}
\log\left(\eta(q^k)^{-1}\sum^{\infty}_{n=-\infty}q^{kn^2+hn}\right)=-\sum^{\infty}_{n=1}f_{k,h}(n)q^{2n}+\sum^{\infty}_{n=1}f_{k,h}(n)q^n.
\end{equation}
\\
\textbf{Proof.}\\
For $|q|<1$, $k>0$, the Jacobi triple product identity ([17] Theorem 13 chapter 8 page 169) becomes
\begin{equation}
\sum^{\infty}_{n=-\infty}q^{kn^2+hn}=\prod^{\infty}_{n=0}\left(1-q^{2kn+2k}\right)\left(1+q^{2kn+k+h}\right)\left(1+q^{2kn+k-h}\right).
\end{equation}
Hence we get
\begin{equation}
\sum^{\infty}_{n=-\infty}q^{kn^2+hn}=\eta(q^k)\prod^{\infty}_{n=1}(1+q^n)^{\chi_{k,h}(n)},
\end{equation}
when $k>|h|>0$, $k\in\bf N\rm$, $h\in\bf Z\rm$.\\
Using the fact that for $x>0$ holds the relation (analogous to (3)):
\begin{equation}
\exp\left(f\left(e^{-x}\right)-f\left(e^{-2x}\right)\right)=\prod^{\infty}_{n=1}\left(1+e^{-n x}\right)^{\frac{1}{n}\sum_{d|n}\frac{f^{(d)}(0)}{\Gamma(d)}\mu\left(\frac{n}{d}\right)}.
\end{equation}
Hence we get (126).\\

Also we have the next general result:\\
\\
\textbf{Theorem 16.}\\
If $|q|<1$ and $k,h$ integers such $k>|h|>0$, then 
\begin{equation}
\sum^{\infty}_{n=-\infty}q^{kn^2+hn}=\exp\left(-\sum^{\infty}_{n=1}G_{k,h}(n)q^n\right),
\end{equation}
where
\begin{equation}
G_{k,h}(n):=f_{k,h}\left(\frac{n}{2}\right)-f_{k,h}(n)+\frac{k}{n}\sigma_1\left(\frac{n}{k}\right)
\end{equation}
and $\sigma_1(n)=\sum_{d|n}d$.\\
(In relation (131) we have assumed that only integer points count, i.e. $y_{k,h}\left(\frac{n}{2}\right)=0$, when $n$ is odd and takes its regular value otherwise (the same thing holds and with $\sigma_1\left(\frac{n}{k}\right)$)).\\
\\
\textbf{Proof.}\\
From 
\begin{equation}
\log\left(\eta\left(q^k\right)\right)=-\sum^{\infty}_{n=1}\frac{q^{nk}}{n}\sum_{d|n}d\textrm{, }|q|<1
\end{equation} 
and relation (126) we have
\begin{equation}
\log\left(\sum^{\infty}_{n=-\infty}q^{kn^2+hn}\right)=-\sum^{\infty}_{n=1}f_{k,h}(n)q^{2n}+\sum^{\infty}_{n=1}f_{k,h}(n)q^n-\sum^{\infty}_{n=1}\frac{q^{nk}}{n}\sum_{d|n}d,
\end{equation}
and the result follows.\\
\\
\textbf{Theorem 18.}\\
If $k_l,h_l$ are integers such $k_l>|h_l|>0$, then
the number of representations of a positive integer $n$ into form
\begin{equation}
n=\sum^{N}_{l=1}\left(k_lX_l^2+h_lX_l\right)
\end{equation}  
is
\begin{equation}
r_{\{k_1,k_2,\ldots,k_N;h_1,h_2,\ldots,h_N\}}(n)=T\left(\sum^{N}_{l=1}G_{k_l,h_l}(n)\right).
\end{equation}
\\
\textbf{Examples.}\\
\textbf{i)} The Diophantine equation
\begin{equation}
2X_1^2+3X_2^2-X_1+2X_2=3
\end{equation}
have $r(3)=r_{\{2,3;-1,2\}}(3)$ solutions. For to find the value of $r(3)$ we use the relations (70),(127),(132),(136):\\
For $k=2$ and $h=-1$, we have $f_{2,-1}(1)=1$, $f_{2,-1}(2)=3/2$, $f_{2,-1}(3)=4/3$.\\
For $k=3$ and $h=2$ we have $f_{3,2}(1)=1$, $f_{3,2}(2)=1/2$, $f_{3,2}(3)=4/3$.\\ 
Hence\\  
$G_{2,-1}(1)=0-1+0=-1$, $G_{2,-1}(2)=1-3/2+1=1/2$, $G_{2,-1}(3)=0-4/3+0=-4/3$. 
Also\\ 
$G_{3,2}(1)=-1$, $G_{3,2}(2)=1/2$, $G_{3,2}(3)=-1/3$.\\
Hence 
$G_{2,-1}(1)+G_{3,2}(1)=-2$, $G_{2,-1}(2)+G_{3,2}(2)=1$, $G_{2,-1}(3)+G_{3,2}(3)=-5/3$
and
$b_0=1$, $b_1=2$, $b_2=1$, $b_3=1$.
Hence equation (137) has 1 solution in the set of integers.\\
Note also that equations
\begin{equation}
2X_1^2+3X_2^2-X_1+2X_2=2
\end{equation} 
\begin{equation}
2X_1^2+3X_2^2-X_1+2X_2=1
\end{equation} 
have exactly 1 and 2 solutions respectively.\\
\textbf{ii)} 
The equation
\begin{equation}
3X_1^2+5X_2^2+6X_3^2-X_1+2X_2-2X_3=20
\end{equation}
have 3 solutions in the set of integers.\\
\\

\[
\]

\centerline{\bf References}\vskip .1in

[1] T. Apostol. 'Introduction to Analytic Number Theory'. Springer Verlag, New York, Berlin, Heidelberg, Tokyo, 1974.

[2] J.V Armitage, W.F Eberlein. 'Elliptic Functions'. 
Cambridge University Press. 2006.

[3] Nikos Bagis. 'Evaluations of Derivatives of Jacobi Theta Functions in the Origin'. arXiv:1105.6279v1 [math.GM]. 2011. 

[4] Nikos Bagis. 'Generalizations of Ramanujan's Continued Fractions'. arXiv: 1107.2393v2 [math.GM]. 2012.

[5] J.M. Borwein and P.B. Borwein. 'Pi and the AGM: A Study in Analytic Number Theory and Computational Complexity', Wiley, New York, 1987.

[6] E.T. Whittaker and G.N. Watson. 'A course on Modern Analysis'. Cambridge U.P. 1927.

[7] B.C. Berndt. 'Ramanujan`s Notebooks Part II'. Springer-Verlag, New York. 1989.

[8] M.L. Glasser and N.D. Bagis. 'Some applications of the Poisson summation formula'. arXiv. 2008.

[9] Bruce C. Berndt. 'Ramanujan`s Notebooks Part III'. Springer-Verlag, New York. 1991.

[10] B.C. Berndt. 'Ramanujan's Notebooks Part V'. Springer Verlag, New York. 1998.

[11] N.D. Bagis. 'A General Method for Constructing Ramanujan-Type Formals for Powers of $1/\pi$'. The Mathematica Journal. Vol. 15. 2013

[12] D. Broadhurst. 'Solutions by radicals at Singular Values $k_N$ from New Class Invariants for $N\equiv3\;\; mod\;\; 8$'. arXiv:0807.2976 [math-ph]. 2008

[13] S. Ramanujan. 'On certain arithmetical functions'. Trans. Cambridge Philos. Soc., Vol. 22, 159-184. 1916 

[14] M. Abramowitz and I.A. Stegun. 'Handbook of Mathematical Functions'. Dover Publications, New York. 1972

[15] N.D. Bagis. 'Some New Results on Sums of Primes'. Mathematical Notes, Vol. 90, No. 1, pp 10-19. 2011

[16] N.D. Bagis, M.L. Glasser. 'Conjectures on the Evaluation of Certain Functions with Algebraic Properties'. Journal of Number Theory. 155, 63-84. 2015

[17] G.E. Andrews. 'Number Theory'. Dover Publications, Inc. New York. 1994

[18] Carlos J. Moreno, Samuel S. Wagstaff, Jr. 'Sums of Squares of Integers'. Chapman and Hall/CRC,. Taylor and Francis Group, Boca Raton, London, New York. 2006  

[19] J.M. Borwein, M.L. Glasser, R.C. McPhedran, J.G. Wan, I.J. Zucker. 'Lattice Sums Then and Now'. Cambridge University Press. 2013

[20] A. Raouf Chouikha. 'Complete Monotonicity of classical theta functions and applications'. arXiv: 1409.1498v1 [math.CA] 4 Sep 2014 

[21]: G.E.Andrews, Amer. Math. Monthly, 86, 89-108,(1979).

\end{document}